\documentclass[12pt]{amsart}
\newcommand{\R}{\mathbb R}

\usepackage{pdfsync}
\usepackage[active]{srcltx}
\usepackage[english]{babel}
\usepackage{amscd}
\usepackage{amssymb}
\usepackage{amsthm}
\usepackage{amsmath}
\usepackage[latin2]{inputenc}
\usepackage{t1enc}
\usepackage{graphicx}
\usepackage{comment}
\usepackage{enumerate}
\usepackage{color}
\usepackage{hyperref}
\usepackage{wrapfig}

\usepackage{psfrag}
\usepackage{epsfig}

\usepackage{dsfont}
\usepackage{enumerate}
\usepackage{bbm}
\usepackage{soul}

\definecolor{trp}{rgb}{1,1,1}

\definecolor{red}{rgb}{1,0,.2}

\definecolor{blue}{rgb}{0,0,1}

\definecolor{rgrey}{rgb}{.8,0.4,.4}  % faint
\definecolor{grey}{rgb}{.13,.13,.13}  % almost black

\definecolor{green}{rgb}{0.0,0.4,0.2}

\newcommand*{\ind}{\mathbf{1}}

\newtheorem{theorem}{Theorem}
\theoremstyle{plain}

\newtheorem{definition}[theorem]{Definition}
\newtheorem{example}[theorem]{Example}

\newtheorem{lemma}[theorem]{Lemma}

\newtheorem{proposition}[theorem]{Proposition}
\newtheorem{remark}[theorem]{Remark}

\numberwithin{equation}{section}

\usepackage[usenames,dvipsnames,svgnames]{xcolor}

\begin{document}

\parindent0pt

\title[Fractal percolations]
{The geometry of fractal percolation,}

\author{Micha\l\ Rams}
\address{Micha\l\ Rams, Institute of Mathematics, Polish Academy of Sciences, ul. \'Sniadeckich 8, 00-956 Warsaw, Poland
\tt{rams@impan.gov.pl}}

\author{K\'{a}roly Simon}
\address{K\'{a}roly Simon, Institute of Mathematics, Technical
University of Budapest, H-1529 B.O.box 91, Hungary
\tt{simonk@math.bme.hu}}

 \thanks{
\indent Rams was partially supported by the MNiSW grant N201 607640 (Poland).
 The research of Simon was supported by OTKA Foundation
\# K 104745}

\begin{abstract}
A well studied family of random fractals called fractal percolation is discussed. We focus on the projections of fractal percolation  on the plane. Our goal is to present stronger versions of the classical Marstrand theorem, valid for almost every realization of fractal percolation. The extensions go in three directions:
\begin{itemize}
\item the statements work for all directions, not almost all,
\item the statements are true for more general projections, for example radial projections onto a circle,
\item in the case $\dim_H >1$, each projection has not only positive Lebesgue measure but also has nonempty interior.
\end{itemize}

%The coordinate axes projections and the $45^{\circ}$ projection (which is related to the arithmetic sum) have been extensively studied (see \cite{Dekking2009},\cite{Dekking1988},\cite{Dekking1990},\cite{Dekking2008},
%\cite{Falconer1986},\cite{Falconer1992},\cite{Falconer1994}).
%  After a short review of these  results, we show that at least in the homogeneous case, for almost all realizations  one of the two things happens:
 % either the dimension of $E$ is smaller than one and then dimension  of $E$ is preserved by \textbf{all} linear \textbf{projections}. Otherwise, \textbf{all} linear \textbf{projections} of $E$ contain some intervals.  This gives much more precise information about the structure of projections of fractal percolation than Marstrand Theorem (Theorem 1).
%We point out that the assertions above hold  for some important non-linear projections like radial and co-radial projections (see Figure \ref{n4}).

 \end{abstract}
%\date{\today}

\maketitle

\medskip

\section{introduction}

To model turbulence, Mandelbrot \cite{Mandelbrot1974,Mandelbrot1983} introduced a
 statistically self-similar family of random Cantor sets.
 Since that time this family has got at least three  names in the literature:  fractal percolation, Mandelbrot percolation and canonical curdling, among which we will use the first one.

In 1996 Lincoln Chayes \cite{Chayes1996} published an excellent survey giving an account about the most important results known in that time. His survey focused on the percolation related properties while we place emphasis on the geometric measure theoretical properties (projections and slices) of fractal percolation sets.

\begin{figure}
  % Requires \usepackage{graphicx}
  \includegraphics[width=8cm]{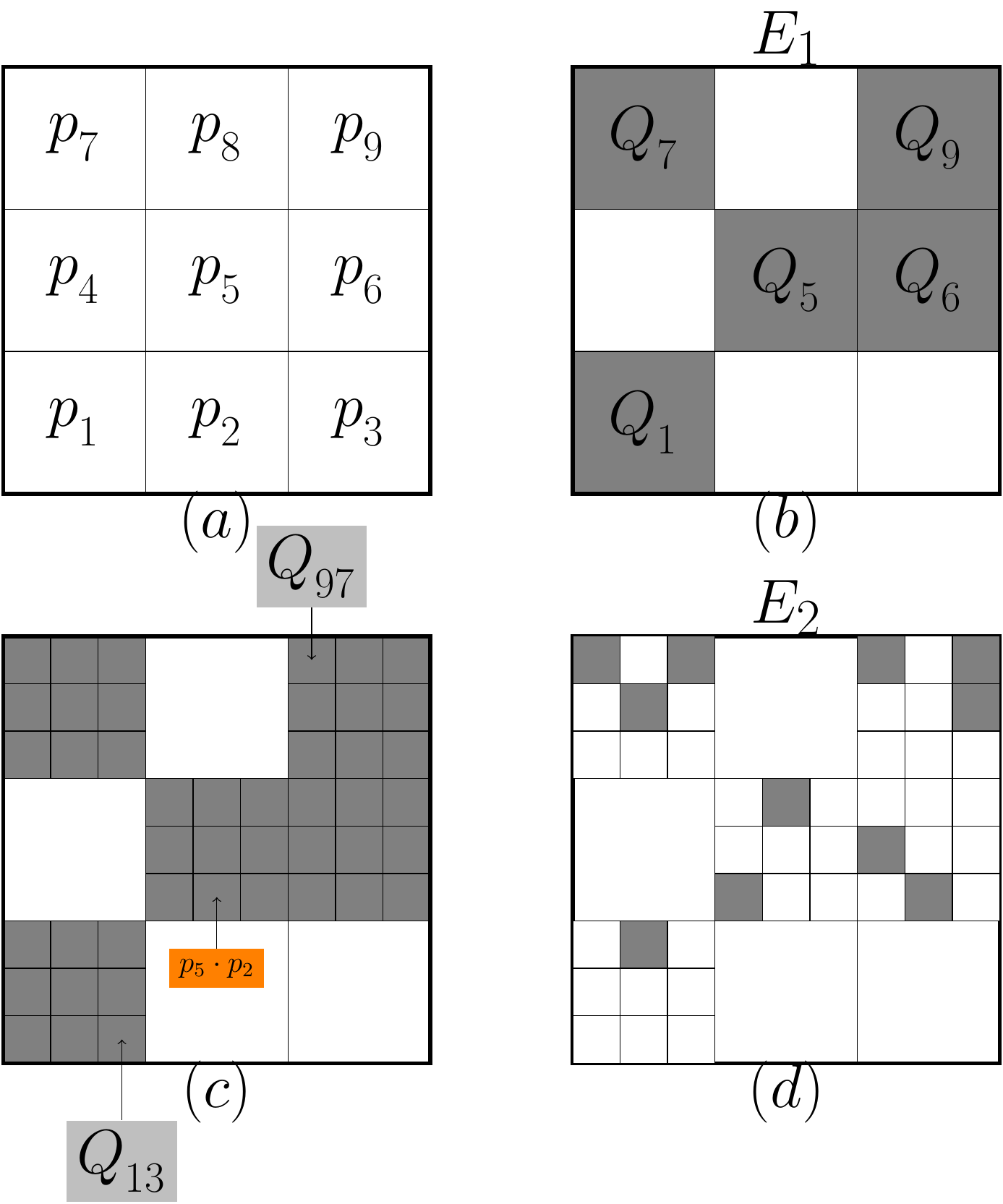}\\
  \caption{The first two steps of the construction. $\mathbb{P}\left(Q_{52}\mbox{ retained }\right)=p_5\cdot p_2$. For this realization $\mathcal{E}_1=\left\{1,5,6,7,9\right\}$, $\mathcal{E}_2=\left\{17,51,58,
  62,64,75,77,79,96,97,99
  \right\}$}.\label{1}
\end{figure}

About the projections of a general Borel set the celebrated Marstrand Theorem gives the following information:

\begin{theorem}[\cite{Marstrand1954}]
Let $E\subset \mathbb{R}^2$ be a Borel set. %We distinguish the cases when $\dim_{\rm H}(E)$ is bigger and when it is smaller than one:
\begin{itemize}
  \item If $\dim_{\rm H}(E)<1$ then for Lebesgue almost all $\theta $ $\dim_{\rm H}(\mathrm{proj}_\theta (E))=\dim_{\rm H}(E)$.
  \item If $\dim_{\rm H}(E)>1$ then for Lebesgue almost all $\theta $ we have $\mathcal{L}{\rm eb}(\mathrm{proj}_\theta (E))>0$.
\end{itemize}
where $\mathrm{proj}_\theta $ is the orthogonal projection in direction $\theta $.
\end{theorem}
In this paper we review some recent results which give more precise information in the special case of the projections of fractal percolation Cantor sets.

\section{The construction and its immediate consequences}

The construction consists of the infinite iteration of two steps. We start from the unit cube in $\mathbb{R}^d$.
\begin{itemize}
  \item All cubes we have after the $n$-th iteration  of the process (they will be called level $n$ cubes) we subdivide into smaller cubes of equal size,
  \item Among them some are retained and some are discarded. Retaining or discarding of different cubes are independent random events. The cubes that were retained are the level $n+1$ cubes.
\end{itemize}
Those points that have never been discarded form the fractal percolation set.

Please note that in literature the term fractal percolation is often used to denote object which we call homogeneous fractal percolation. That is, the fractal percolation for which all squares have equal probabilities of being retained.

\subsection{An informal description of Fractal Percolation}

We fix integer $M\geq 2$.
We partition the unit cube  $Q\subset \mathbb{R}^d$ into $M^d$ congruent cubes of side length $M^{-1}$ and we assign a probability to each of the cubes in this partition (Figure \ref{1} (a)).
We retain each of the cubes of this partition with the corresponding probability independently and discard it with one minus the corresponding probability.
The union of the retained squares is the first approximation of the random set to be constructed ( Figure \ref{1} (b)). We obtain the second approximation by repeating this process independently of everything in each of the retained squares ( Figure \ref{1} (c) and (d)). We continue this process at infinitum.

The object of our investigation is
the collection of those points which have not been discarded. It will be called \textbf{fractal percolation set} and denoted by ${E=(d,M,\mathbf{p})}$, where $\mathbf{p}$ is the chosen vector of the probabilities $\{p_i\}$.
In the special case when all $p_i$ are equal we obtain the \textbf{homogeneous fractal percolation set} which is denoted by $E^h=E^h(d,M,p)$.

\subsection{Fractal percolation set in more details}
For simplicity we give the construction on the plane but the definition works with obvious modifications in $\mathbb{R}^d$ for all $d\geq 1$. Besides the dimension of the ambient space the two other parameters of the construction are: the natural number $M\geq 2$ and
a vector of probabilities $\mathbf{p}\in [0,1]^{M^2}$ (note: not a probabilistic vector).
To shorten the notation we write $\mathcal{I}$ for the set of indices of $\mathbf{p}$:
$$
\mathcal{I}:=\left\{1,\dots ,M^2\right\}
$$
The statistically self-similar random  set which is the object of our study is defined as
\begin{equation}\label{4}
  E:=\bigcap _{n=1}^{\infty }E_n,
\end{equation}
where $E_n$ is the $n$-th approximation of $E$. The inductive definition of $E_n$ will occupy the rest of this subsection. Actually $E_n$ is the union of a random collection of level $n$ squares. First we define the level $n$ squares and then we introduce the random rule with which  those level $n$ squares are selected whose union form $E_n$.

\subsubsection{The process of subdivision}
We divide the unit square $Q=\left[0,1\right]^2$ into $M^2$ congruent squares $Q_1,\dots ,Q_{M^2}$ of size $M^{-1}$ numbered according to lexicographical order (or any other order). These squares are the level one $M$-adic squares. Let
$$
\mathcal{N}_1:=\left\{x_i\right\}_{i\in \mathcal{I}}
$$
be the set of midpoints of the level one squares. For each midpoint $x_i$ we define the homothetic map  $\varphi _i:Q\to Q_i$:
$$
\varphi _i(y):=x_i+M^{-1}\cdot \left(y-\left(\frac{1}{2},\frac{1}{2}\right)\right).
$$
For every $\mathbf{i}\in \mathcal{I}^n$, $\mathbf{i}=(i_1,\dots ,i_n)$ we write
$$
x_{\mathbf{i}}:=\varphi _{\mathbf{i}}\left(\frac{1}{2},\frac{1}{2}\right).
$$
and we define the map
$$
\varphi _{\mathbf{i}}(y):=x_{\mathbf{i}}+M^{-n}\cdot \left(y-\left(\frac{1}{2},\frac{1}{2}\right)\right).
$$
To simplify the notation, we will not distinguish the set of the centers of level $n$ squares
$$
\mathcal{N}_n:=\left\{\varphi _{\mathbf{i}}\left(\frac{1}{2},\frac{1}{2}\right):\mathbf{i}\in \mathcal{I}^n\right\}
$$
and the family of level $n$-squares:
\begin{equation}\label{3}
  \left\{Q_{\mathbf{i}}:=\varphi _{\mathbf{i}}(Q): \mathbf{i}\in \mathcal{I}^n\right\}.
\end{equation}
\begin{figure}\label{n3}
  % Requires \usepackage{graphicx}
  \includegraphics[width=12cm]{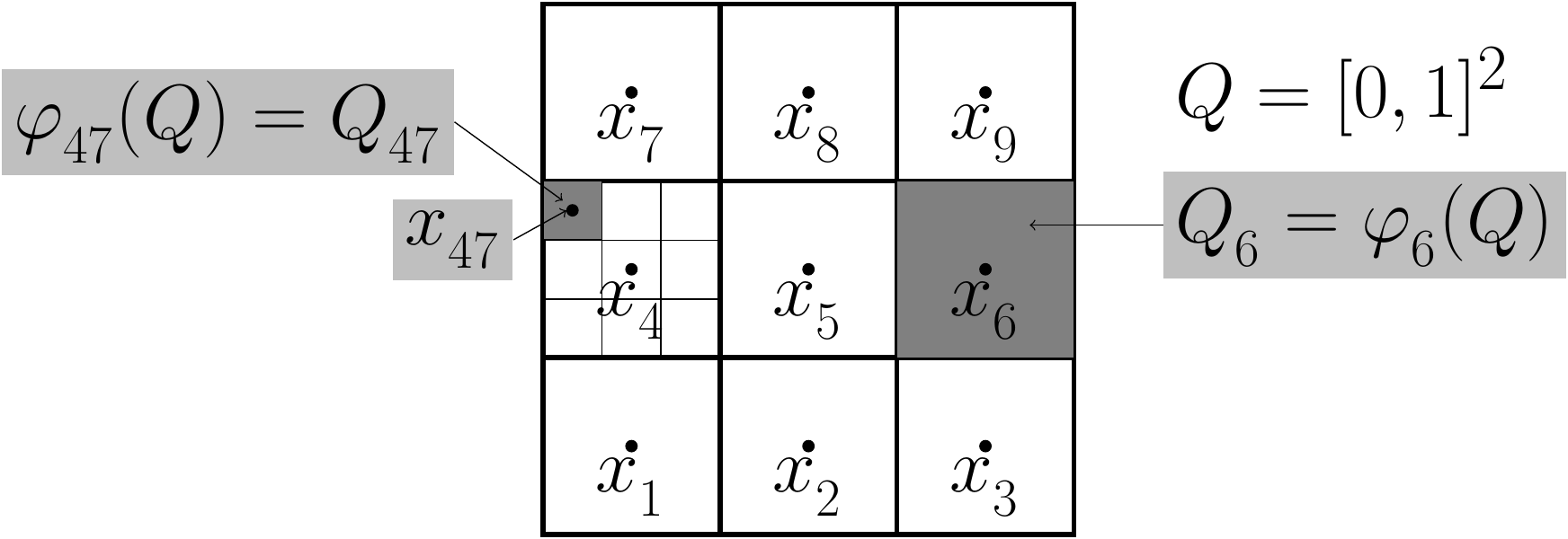}\\
 \caption{Definition of level $n$ squares}\label{2}
\end{figure}
\subsubsection{The process of retention}

%For each $i\in \mathcal{I}$ we have a coin that we call the $i$-th coin. It lands head with probability $p_i$ (which was the $i$-th component of the given vector $\mathbf{p}$). Now for each $i\in \mathcal{I}$ we flip the $i$-th coin independently. Those level one squares  $Q_i$ where the coin flipping resulted in head are retained, the others are discarded.

The square $Q=Q_{\emptyset}$ is retained.
For any $\mathbf{i}\in \mathcal{I}^n$ for which the square $Q_{\mathbf{i}}$ is retained and for each $j\in \mathcal{I}$, the square $Q_{\mathbf{i} j}$ is retained with probability $p_j$. The events '$Q_{\mathbf{i} j}$ is retained' and '$Q_{\mathbf{i'} j'}$ is retained' are independent whenever $\mathbf{i} \neq \mathbf{i'}$ or $j\neq j'$.

We define $E_1$ as the union of retained squares $Q_i, i\in \mathcal{I}$. Similarly, $E_n$ is the union of retained squares $Q_\mathbf{i}, \mathbf{i}\in \mathcal{I}^n$. We write
$$
\mathcal{E}_n:=\left\{\mathbf{i}\in \mathcal{I}^n:Q_\mathbf{i}\mbox{ retained }\right\}.
$$

\subsection{The corresponding probability space and statistical self-similarity}
The  probability space corresponding to this random construction is best described by M.
Dekking  \cite{Dekking2009}. For the convenience of the reader we repeat it here.
Let  $\mathcal{T}$ be the $M^d$ array tree that is
$$
\mathcal{T}:=\bigcup _{n=0}^{\infty }
\mathcal{I}^n,
$$
where $\mathcal{I}^0:=\emptyset $ is the root of three . Let $\Omega :=\left\{0,1\right\}^ {\mathcal{T}}$ that is $\Omega $ is the set of labeled trees where we label every node of $\mathcal{T}$ by $0$ or $1$. The probability measure $\mathbb{P}_{\mathbf{p}}$ on $\Omega $
is define in such a way that the family of labels  $X_{\mathbf{i}}\in \left\{0,1\right\}$ of nodes $\mathbf{i}\in \mathcal{T}$ satisfy:
\begin{itemize}
  \item $\mathbb{P}_{\mathbf{p}}(X_\emptyset =1)=1$
  \item $\mathbb{P}_{\mathbf{p}}(X_{i_1,\dots ,i_n})=p_{i_n}$
  \item  $\left\{X_{\mathbf{i}}\right\}_{\mathbf{i}\in \mathcal{T}}$ are independent.
\end{itemize}
Following \cite{Dekking2009} we define the survival set of level $n$ by
$$
S_n:=\left\{\mathbf{i}\in \mathcal{I}^n:X_{i_1\dots ,i_k}=1,\ \forall 1\leq k\leq n\right\}.
$$
Then
$$
E_n=\bigcup _{\mathbf{i}\in S_n} Q_{\mathbf{i}},\quad E=\bigcap _{n=1}^{\infty }E_n.
$$
It follows from the construction that generalized fractal percolation set is statistically self-similar and the number of retained cubes form a branching process:
\begin{lemma}
\begin{description}\label{5}\
  \item[(a)] $\left\{\#\mathcal{E}_n\right\}$ is a branching process with average number of offsprings $\sum\limits_{i\in \mathcal{I}}p_i$. In particular if $p_i\equiv p$ then the
   offspring distribution is $\texttt{Binomial}( M^d,p)$.
  \item[(b)] For every $n\geq 1$ and $\mathbf{i}\in \mathcal{E}_n$ the rescaled copy $\varphi _{\mathbf{i}}^{-1}(E\cap Q_{\mathbf{i}})$ has the same distribution as $E$ itself.
  \item[(c)] The sets $\left\{E\cap Q_{\mathbf{i}}\right\}_{\mathbf{i}\in \mathcal{E}_n}$ are independent.
\end{description}
\end{lemma}
Using this it is not hard to prove that
\begin{equation}\label{6}
  E\ne\emptyset \mbox{ implies that }
\dim_{\rm H}(E)=\dim_{\rm B}(E)=\frac{\log \sum\limits_{i\in \mathcal{I}}p_i}{\log M}\mbox{ a.s. }
\end{equation}
This was proved by Kahane and Peyriere
\cite{Kahane1976}, Hawkes \cite{Hawkes1981}, Falconer \cite{Falconer1986}, Mauldin and Williams \cite{Mauldin1986} independently. A canonical example of the inhomogeneous fractal percolation set is:
\begin{example}[Random Sierpi\'nski Carpet ]
Let $SC_p:=E(2,3,\mathbf{p})$, where using the notation of Figure \ref{1} (c):
$$
p_5=0 \mbox{ and for }i\in \left\{1,\dots ,9\right\}\setminus\left\{5\right\}:\
p_i=p.
$$
\end{example}

\section{Percolation and projection to coordinate axes}
In this section we work on the plane so $Q=[0,1]^2$. The connectivity properties of $E^h(2,M,p)$ for an arbitrary $M\geq 2$  was first investigated by Chayes, Chayes and Durrett \cite{Chayes1988}. Dekking and Meester \cite{Dekking1990} gave a simpler proof and extended the scope of the theorem for some inhomogeneous fractal percolation sets like the random Sierpi\'nski carpet $SC_p$. Here we   summarize briefly some of the most interesting results of this area. For a  much more detailed account see by L. Chayce \cite{Chayes1996}.

We say that $E$ \textbf{percolates} if $E$ contains a connected set which intersects both the left and the right sides of  $Q$. If $E$ percolates then $E$ has a large connected component.

\subsection{The homogeneous case}
The following very important result was proved by Chayce,Chayce, Durrett.
\begin{theorem}[\cite{Chayes1988}] Fix an arbitrary $M\geq 2$. Then there is a critical probability $\frac{1}{M}<p_c<1$ such that
\begin{enumerate}
  \item If $p<p_c$ then $E^h(2,M,p)$ is a random dust that is totally disconnected almost surely.
  \item If $p\geq p_c$ then $E^h(2,M,p)$ percolates with positive probability. This implies that $E^h(2,M,p)$ is not totally disconnected almost surely.
\end{enumerate}
\end{theorem}
This shows a remarkable difference in between the fractal percolation and the usual percolation: in the latter case, the probability of percolation at critical parameter $p=p_c$ is 0.

\subsection{The inhomogeneous case}
Using some earlier works of Dekking and Grimmett \cite{Dekking1988}, the results above were extended by Dekking and Meester \cite{Dekking1990}. They proved that by changing  the components of $\mathbf{p}$ the inhomogeneous fractal percolation set $E(2,M,\mathbf{p})$ can go through the  six stages below. Here the projection to the $x$-axis is denoted by $\mathrm{proj}_x$. That is $\mathrm{proj}_x(a,b)=a$.

\medskip

\textbf{The DM stages of} $\mathbf{E(2,M,\mathbf{p})}$:
\begin{description}
  \item[I]   $E=\emptyset $ almost surely.
  \item[II] $\mathbb{P}\left(E\ne\emptyset \right)>0$ but $\dim_{\rm H}\left(\mathrm{proj}_xE\right)=\dim_{\rm H}\left(E\right)$ almost surely.
  \item[III] $\dim_{\rm H}\left(\mathrm{proj}_xE\right)<\dim_{\rm H}\left(E\right)$ if $E\ne\emptyset $ but $\mathcal{L}{\rm eb}\left(\mathrm{proj}_xE\right)=0$ almost surely.
  \item[IV] $0<\mathcal{L}{\rm eb}\left(\mathrm{proj}_xE\right)<1$ almost surely.
  \item[V] $\mathcal{L}{\rm eb}\left(\mathrm{proj}_xE\right)=1$ holds with positive probability but $E$ does not percolate almost surely.
  \item[VI] $E$  percolates with positive probability.
\end{description}
It was proved in \cite{Dekking1990} that the random Sierpi\'nski Carpet $SC_p$ goes through all of these stages as we increase the value of $p$. The following theorem gives the precise answer when exactly a system appears in stages I,II,III.
\begin{theorem}[\cite{Dekking1988}, \cite{Falconer1986}]   Let $m_r$ be the sum of the probabilities in the $r$-th column, that is the expected number of squares in column $r$. Then
\begin{enumerate}
  \item $E=\emptyset $ almost surely iff  $\sum\limits_{i=1}^{M^2}p_i\leq 1$. Except when $\exists i$ such that $p_i=1$ and $p_j=0$ for all $i\ne j$. In this case $E$ is a singleton.
  \item $\dim_{\rm H}(\mathrm{proj}_x(E))=\dim_{\rm H}(E)$ holds almost surely, iff $\sum\limits_{r=1}^{M}m_r \log m_r\leq 0$.
  \item $\mathcal{L}{\rm eb}(\mathrm{proj}_xE)=0$ holds almost surely iff $\sum\limits_{r=1}^{M}\log m_r\leq 0$.
\end{enumerate}
\end{theorem}
This result was strengthened by  Falconer and Grimmett:
\begin{theorem}[\cite{Falconer1992},\cite{Falconer1994}]\label{n9}
Assume that $m:=\min\left\{m_r\right\}>1$. Then $\mathrm{proj}_x(E)$ contains an interval almost surely, conditioned on non-extinction.
\end{theorem}

We will present the proof in the fifth section.

\subsection{The DM stages for the homogeneous case}

For the homogeneous case $m_r=M\cdot p$ Hence we obtain that almost surely:
\begin{itemize}
  \item If $0<p\leq \frac{1}{M^2}$ then $E=\emptyset $.
  \item If $\frac{1}{M^2}<p\leq \frac{1}{M}$ then the system is in stage II.
  \item If $\frac{1}{M}<p<p_c$ then the system is in stage V.
\end{itemize}
Stages III and IV do not appear in the homogeneous case.

\section{The arithmetic sum/difference of two fractal percolations}
There is a very nice and more detailed survey of this field due to M. Dekking \cite{Dekking2009}.
In the previous section we studied the connectivity properties and the $90^{\circ}$ projections of random Cantor sets. In this section we consider sets which are products of inhomogeneous fractal percolation sets and we take their $45^{\circ} $,  ($-45^{\circ} $) projections in order to study the arithmetic difference (arithmetic sum) respectively of independent copies of $E(1,M,\mathbf{p})$.

\subsection{The arithmetic sum and its visualization}
\begin{figure}[h!]
  % Requires \usepackage{graphicx}
  \includegraphics[width=6cm]{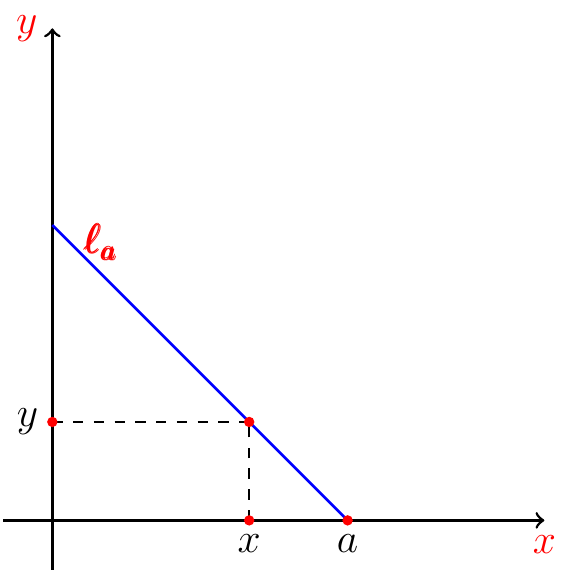}\\
\caption{Algebraic sum as $-45^{\circ}$ projection: $a=x+y=\mathrm{proj}_{-45^{\circ}}(x,y)$
}\label{7}
\end{figure}
Let $A,B\subset \mathbb{R}$ be arbitrary.  Then the arithmetic sum $A+B:=\left\{a+b:a\in A,b\in B\right\}$ is the $-45^{\circ}$-projection of $A\times B$ to the $x$-axis (this is the direction of the line $\ell_a$ on Figure \ref{7}). Similarly,  we can visualize the arithmetic difference by taking the projection of the product set with the line of $+45^{\circ}$ angle.

The motivation for studying the arithmetic difference (or sum) of random Cantor sets comes from a conjecture of Palis which states that typically (in a natural sense which depends on the actual setup), the  arithmetic difference of two dynamically defined Cantors is either small in the sense that it has Lebesgue measure zero or big in the sense that it contains some intervals, but at least typically, it does not occur that the arithmetic difference set is a set of positive Lebesgue measure with empty interior. This conjecture does not hold for the algebraic difference of inhomogeneous fractal percolation sets, but it holds in the homogeneous case. The way to prove this is via the $45^{\circ}$-projections of $E(1,M,\mathbf{p})\times E(1,M,\mathbf{p})$.

\subsection{The product of two one dimensional fractal percolation versus a two dimensional fractal percolation.}\label{8}

We explain this relation in the case when $M=3$.
Assume that we are given the inhomogeneous fractal percolations $E(1,3,\mathbf{a})$, and $E(1,3,\mathbf{b})$, where $\mathbf{a}=(a_1,a_2,a_3)$, $\mathbf{b}=(b_1,b_2,b_3)$ are the vectors of probabilities. We define the vector $\mathbf{p}\in [0,1]^9$ as their product $\mathbf{p}=\mathbf{a}\bigotimes \mathbf{b}$ in the natural way which is suggested by  looking  at Figure \ref{1} (a). That is:
$$
p_i:=a_{u}\cdot b_{v}\mbox{ if }
i-1=3*(v-1)+(u-1),\ 1\leq u,v\leq 3.
$$
The reason that $E(2,3,\mathbf{p})$ and $E(1,3,\mathbf{a})\times E(1,3,\mathbf{b})$ are similar is explained in \textbf{(a)} and the essential difference between them is pointed out in \textbf{(b)} below:
\begin{description}
  \item[(a)] Let $\mathbf{i}\in \left\{1,\dots ,9\right\}^n$. Then the probability that $Q_{\mathbf{i}}$ is retained is the same  during the construction of $E(2,3,\mathbf{p})$ and the construction of $E(1,3,\mathbf{a})\times E(1,3,\mathbf{b})$.
  \item[(b)] Let $K$ and $L$ be level $n$ squares for some $n$. Assume that both $K$ and $L$ are retained during the construction of $E(2,3,\mathbf{p})$ and $E(1,3,\mathbf{a})\times E(1,3,\mathbf{b})$. Then
      \begin{itemize}
        \item In the construction of  $E=E(2,3,\mathbf{p})$ the sets $E\cap K$ and $E\cap L$ are independent.
        \item In the construction of $E(1,3,\mathbf{a})\times E(1,3,\mathbf{b})$ the sets $E\cap K$ and $E\cap L$ are independent iff  $\mathrm{proj}_xK\ne \mathrm{proj}_xL$ and  $\mathrm{proj}_yK\ne\mathrm{proj}_yL$ hold.
      \end{itemize}
\end{description}
In dimension $d\geq 2$ the analogy is the same: the probability of the retention of a level $n$ cube is the same for the $d$-dimensional percolation and for the $d$-fold product of the corresponding one dimensional percolations. On the other hand, the future of what ever happens in two distinct retained level $n$ cubes is:
\begin{itemize}
  \item always independent in the $d$-dimensional percolation case,
  \item independent for the  $d$-fold product of the corresponding one dimensional fractal percolations iff the two cubes do not share any common projections to  coordinate axes.b
\end{itemize}

\subsection{The existence of an interval in the arithmetic difference set}\label{n6}

Let $E_1:=E(1,M,\mathbf{p})$ and $E_2:=E(1,M,\mathbf{q})$. We define the cyclic cross correlation coefficients:
\begin{equation}\label{9}
  \gamma _k:=\sum\limits_{i=1}^{M} p_iq_{i-k (\mathrm{mod}\ M)} \mbox{ for } k=1,\dots ,M.
\end{equation}

\begin{theorem}[\cite{Dekking2008}]\label{12}
Assuming that $E_1,E_2\not=\emptyset $, we have
\begin{description}
    \item[(a)] If
    $\forall i=1,\dots ,M: \ \gamma_i>1$
     then almost surely
     $$E_2-E_1 \mbox{ contains an
    interval }.$$
    \item[(b)] If \  $\exists i\in \left\{1,\dots
    ,M\right\}: \   \gamma_i,\gamma_{i+1 \bmod M}<1$ then almost surely
    $$E_2-E_1 \mbox{ does not contain any interval }.$$
\end{description}\end{theorem}

In the homogeneous case and in the case when $M=3$ this gives complete characterization. Otherwise we can change to higher order Cantor sets (collapsing $n\geq 2$ steps of the construction into one) and we can apply the same theorem in that case. The fact that this can be done is not trivial because higher order fractal percolations are correlated. That is the way as the random set develops in one level $n$ square is dependent how it develops in some other squares. Nevertheless, M. Dekking and H. Don proved that this can be done by pointing out that the proof of the theorem above can be carried out for more general, correlated random sets than the inhomogeneous fractal percolations. This more general family includes the higher order fractal percolation sets.

\subsection{The Lebesgue measure of the arithmetic difference set}
Let $E_2,E_2$ be two independent realizations of $E(1,M,\mathbf{p})$. Then
$$
\gamma _k:=\sum\limits_{i=1}^{M} p_ip_{i-k (\mathrm{mod}\ M)} \mbox{ for } k=1,\dots ,M.
$$
Let $\Gamma :=\gamma _1\cdots \gamma _M$.
\begin{theorem}[\cite{Mora2009}]\label{10}
If $\Gamma >1$ then \newline
$\mathcal{L}{\rm eb}(E_2-E_1)>0$.
\end{theorem}

Combined application of Theorems \ref{12} and \ref{10} yields that the Palis conjecture does not hold in the case when
for $M=3$ and $\mathbf{p}=(0.52,0.5,0.72)$.
 Namely, in this case $\gamma _1=1.0388$ and $\gamma _2=\gamma _3=0.941$.
 Let $E_1,E_2$ be two independent copies of $E(1,3,\mathbf{p})$. Then by Theorem \ref{12} there is no interval in $E_1-E_2$ (since there are two consecutive $\gamma $'s that are smaller than one) and by Theorem \ref{10} we have $\mathcal{L}{\rm eb}\left(E_1-E_2\right)>0$ since $\gamma _1\cdot \gamma _2\cdot \gamma _3=1.0272>1$.

\newpage \section{General projections: the opaque case}

In this and in the following sections we study the projections of fractal percolation sets in general directions. In this section we consider the case $\dim_{\rm H}(E)>1$. Under some mild assumption, almost surely projections of $E$ have not only positive Lebesgue measure, as per Marstrand theorem, but also non-empty interior. Furthermore it holds for all and not only almost all directions. Moreover, this remains valid if we replace the orthogonal projection with a much more general family of projections.

%%%%%%%%%%%%%%%%%%%%%%%%%%%%%%%%%%%%%
%%%%%%%%%%%%%%%%%%%%%%%%%%%%%%%%%%%%%

         \begin{figure}[ht!]
  % Requires \usepackage{graphicx}
  \includegraphics[width=8cm]{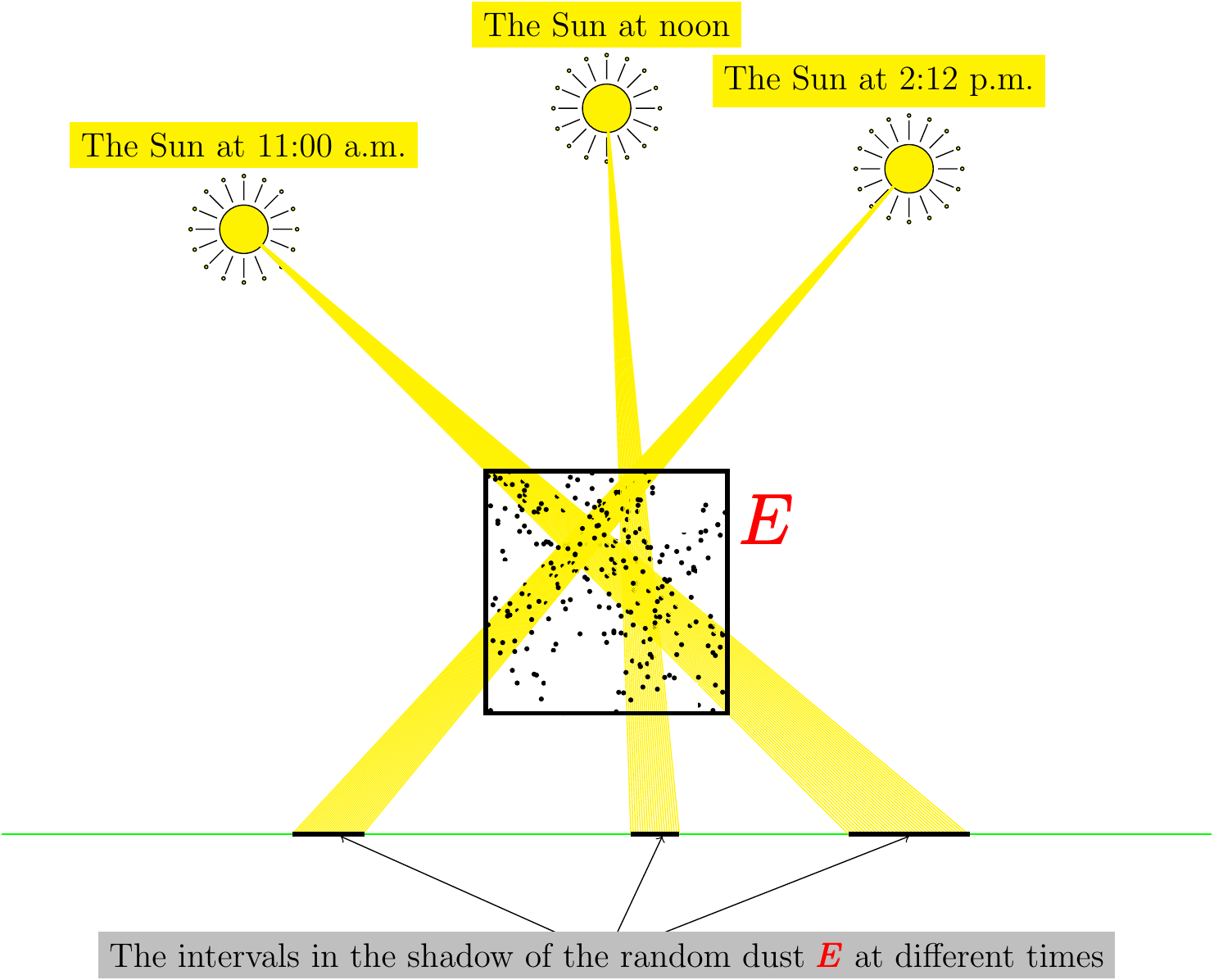}\\
 \caption{}\label{n2}
\end{figure}
One practical application
of our result is shown above (Figure \ref{n2}).
One does not need to rotate such a set to use it as an umbrella.

\begin{figure}\label{n4}
\begin{center}
\includegraphics[width=13.5cm]{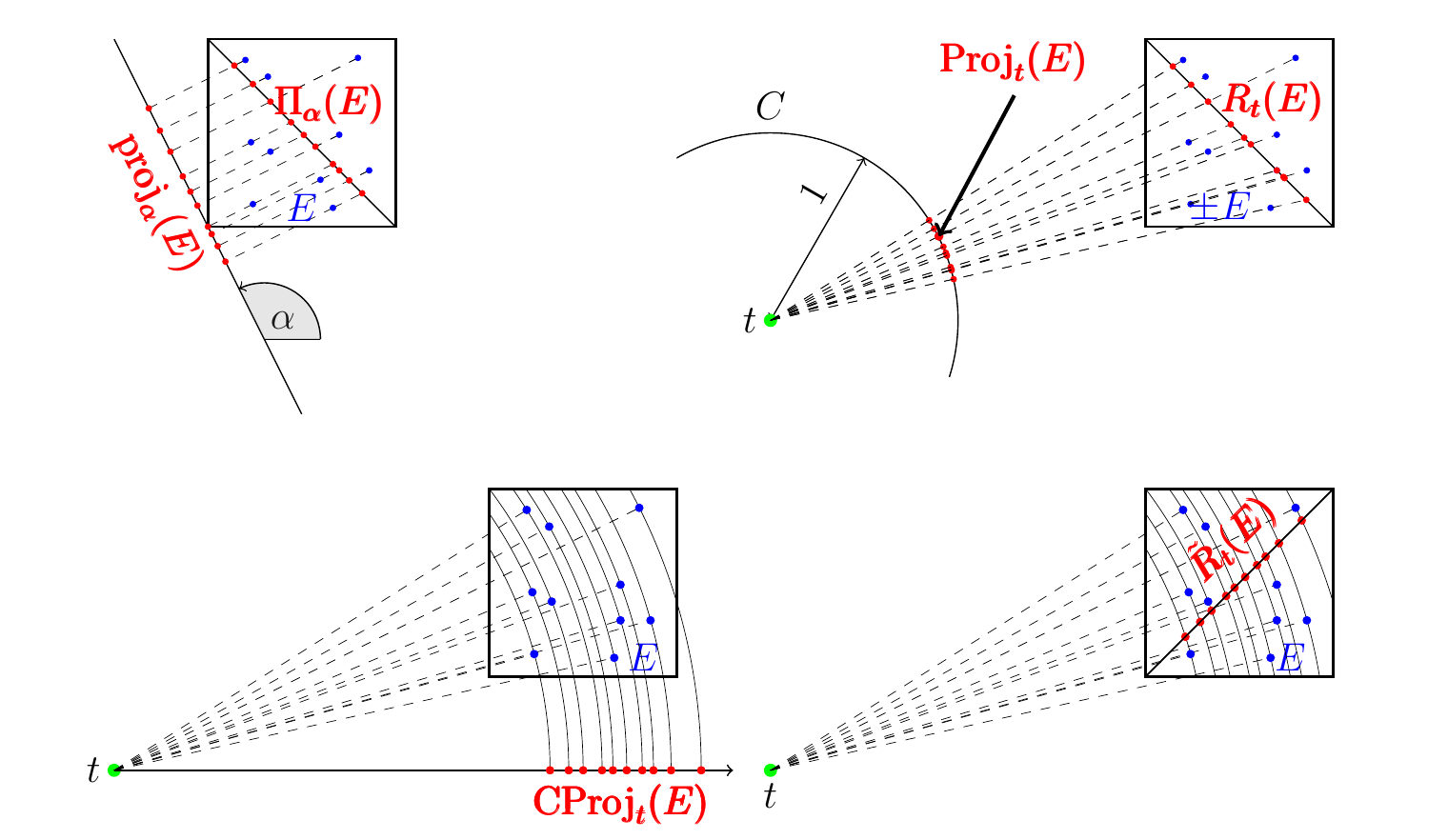}
\caption{The orthogonal $\mathrm{proj}_\alpha $, radial $\mathrm{Proj_t}$, co-radial $\mathrm{CProj}_t$ projections  and the auxiliary projections $\Pi _\alpha $, $R_t$, and $\tilde{R}_t$. }
\label{151}
\end{center}
\end{figure}
We have already studied the horizontal and vertical projections. So we can restrict our attention to the directions  $\alpha \in \mathcal{D}:=\left(0,90^{\circ}\right)$
A condition $A(\alpha )$, $\alpha \in \mathcal{D}$ on the vector of probabilities $\mathbf{p}$ will be defined below.

\begin{theorem}[\cite{Rams}] \label{thm:ortho}
Let $\alpha\in\mathcal{D}$. If $A(\alpha)$ holds and $E$ is nonempty then almost surely ${\rm proj}_\alpha(E)$ contains an interval.
\end{theorem}

 \begin{theorem}[\cite{Rams}] \label{n1}
         If $A(\alpha )$ holds for all $\alpha \in \mathcal{D}$ and $E\ne\emptyset $
         then almost surely all projections ${\rm proj}_\alpha(E)$ contain an interval.
 \end{theorem}

\begin{remark}
The assertions of Theorems \ref{thm:ortho} and \ref{n1} remain valid if we replace $\mathrm{proj}_\alpha $ with more general families of projections, see \cite[Section 6]{Rams}. In particular, radial or co-radial projections (see Figure \ref{151}) are included.
\end{remark}

\begin{example}\label{216}If either
\begin{enumerate}
  \item  Homogeneous case: $p_i= p>M^{-1}$ for all $i$, or
  \item   Generalized random Sierpi\'nski Carpet: $M=3$, $p_{5}=q$, $p_{i}=p$ for $i\ne 5$, and
\[
p>\max\left(\frac 13, \frac {1-q} 2\right)
\]
\end{enumerate}
then Condition A($\alpha$) is satisfied for all $\alpha\in \mathcal{D}$.
Note that (1) is equivalent to $\dim_{\rm H}(E)>1$ almost surely.
\end{example}

\newpage

\subsection{Horizontal and vertical projections}\label{n7}

\begin{wrapfigure}{r}{0.3\textwidth}
  \begin{center}
    \includegraphics[width=3cm]{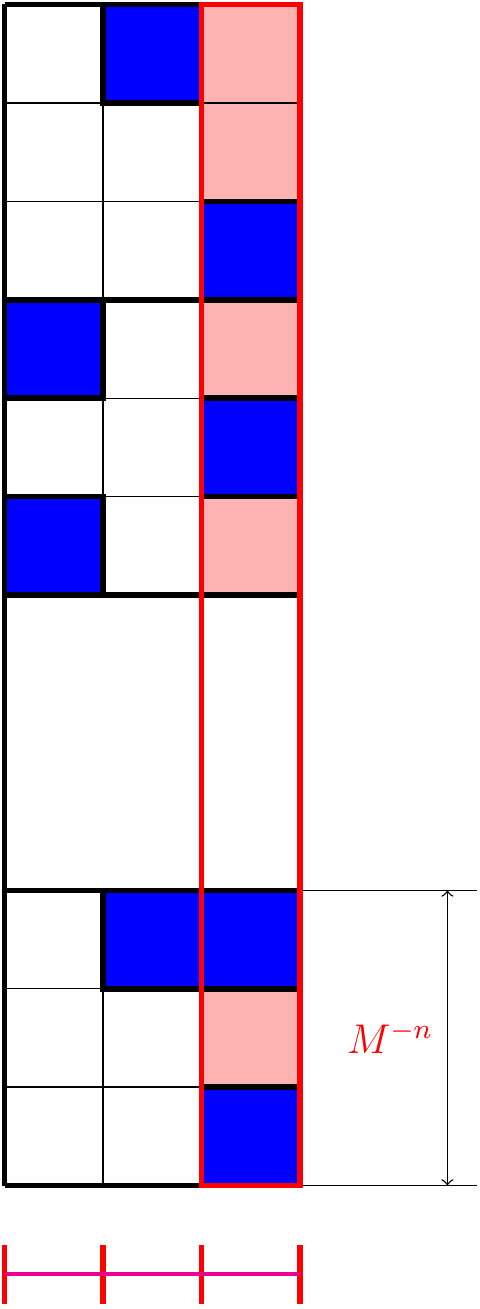}
  \end{center}
 % \caption{fan}
\end{wrapfigure}Let us start by presenting the large deviation argument (by Falconer and Grimmett) working for horizontal and vertical projections. If
$\dim_{\rm H}\Lambda >1$ then from the dimension formula for some $n$ one can find a level $n$ column with (exponentially) many squares.
We prove inductively that in its every $N$-th level sub-column, $N>n$, we typically have exponentially many squares on each level (probability of existence of $N>n$ and an $N$-th level subcolumn which does not have exponentially many squares is super-exponentially small). When we move from level $n$ column to its level $n+1$ subcolumns, each square in the column gives birth to an expected number of $pM>1$ number of level $n+1$ squares in each of the subcolumns. By large deviation theorem there is only a superexponentially small probability that the number of level $n+1$
squares in a subcolumn is smaller than a fixed $\alpha \in(1,pM)$ multiple of the level $n$ squares in the column. By induction, if this exceptional situation does not happen (or happens only finitely many times), for each $N>n$ the number of squares of level $N$ in each subcolumn will be at least of order $\alpha^{N-n}$.

\subsection{Condition A}

Our goal in this subsection is to modify this argument to work in a more complicated situation of projections in general directions. Indeed, contrary to the horizontal/vertical projections case, here it is in general not true that if a line intersects a square of level $n$ then the expected number of squares of level $n+1$ it intersects is greater than 1. It is still true if the line intersects 'central' part of the square, but not if it hits it close to the corners.

Nevertheless, we are able to find a modified version of the argument. We fix $\alpha\in\mathcal{D}$. We are going to consider $\Pi_\alpha$ instead of ${\rm proj}_\alpha$, i.e. we are projecting onto a diagonal $\Delta_\alpha$ of $Q$, see Figure \ref{n4}. For any $\mathbf{i}\in \mathcal{I}^n$ the map $\Pi_\alpha\circ \varphi_{\mathbf{i}}: \Delta_\alpha \to \Delta_\alpha$ is a linear contraction of ratio $M^{-n}$. We will use its inverse: a map $\psi_{\alpha,\mathbf{i} }: \Pi_\alpha(Q_\mathbf{i})\to \Delta_\alpha$. It is a linear expanding map (of ratio $M^n$) and it is onto.

Consider the class of nonnegative real functions on $\Delta_\alpha$, vanishing on the endpoints. There is a natural random inverse Markov operator $G_\alpha$ defined as

\[
G_{\alpha}f(x) = \sum_{i\in \mathcal{E}_1; x\in \Pi_\alpha(Q_i)} f\circ \psi_{\alpha, i}(x).
\]
The corresponding operator on the  $n$-th level is
$$
G^{(n)}_{\alpha}f(x)=
\sum\limits_{\mathbf{i}\in \mathcal{E}_n; x\in \Pi_\alpha(Q_\mathbf{i}} f\circ\psi _{\alpha ,\mathbf{i}}(x).
$$
In particular for any $H\subset \Delta^\alpha $ we have
\begin{equation*}\label{70}
 G_{\alpha}^{(n)}\ind_H(x)=
 \#\left\{\mathbf{i}\in \mathcal{E}_n:
 x\in \Pi _\alpha \left(\varphi _{\mathbf{i}}(H)\right)\right\}.
\end{equation*}
Although $G^{(n)}_\alpha$ should not be thought of as the $n$-th iterate of $G_\alpha$, the expected value of $G^{(n)}_\alpha$
is the $n$-th iterate of the expected value of $G_\alpha$. Namely, let
\begin{equation*}\label{92}
F_\alpha=\mathbb{E}\left[G_{\alpha}\right] \mbox{ and } F_\alpha^n=\mathbb{E}\left[G_{\alpha}^n\right]
\end{equation*}
We then have the formulas

\[
F_\alpha f(x) = \sum_{i\in \mathcal{I}; x\in \Pi_\alpha(Q_{i})} p_{i}\cdot f\circ \psi_{\alpha, i}(x)
\]
and

\begin{equation*}\label{3}
F^n_\alpha f(x)=\sum\limits_{
\mathbf{i}; x\in \Pi_\alpha(Q_\mathbf{i})}
p_{\mathbf{i}}\cdot f\circ \psi _{\alpha ,\mathbf{i}}(x),
\end{equation*}

where

\[
p_{\mathbf{i}} = \prod_{k=1}^n p_{i_k}.
\]
%Hence, $F_\alpha^n$ is indeed the $n$-th iteration of $F_\alpha$ (which explains why we are allowed to use this notation).

\begin{definition} \label{def:conda}
We say the percolation model satisfies {\bf Condition A($\alpha$)} if there exist closed intervals $I_1^\alpha, I_2^\alpha\subset \Delta_\alpha$ and a positive integer $r_\alpha$ such that
\begin{itemize}
\item[i)] $I_1^\alpha\subset {\rm int}I_2^\alpha, I_2^\alpha \subset {\rm int}\Delta_\alpha$,
\item[ii)] $F_\alpha^{r_\alpha}\ind_{I_1^\alpha} \geq 2 \cdot \ind_{I_2^\alpha}$.
\end{itemize}
\end{definition}

It will be convenient to use additional notation. For $x\in\Delta_\alpha$, $\alpha\in \mathcal{D}$, and $I\subset \Delta_\alpha$ we denote

\[
D_n(x,I, \alpha) = \{\mathbf{i} \in \mathcal{I}^n; x\in \Pi_\alpha\circ \varphi_{\mathbf{i}}(I)\}
.
\]
That is, if we write $\ell^\alpha (x)$ for the line segment through $x\in \Delta_\alpha  $ in direction $\alpha $, $D_n(x,I, \alpha)$ is the set of $\mathbf{i}$
for which $\ell^\alpha (x)$ intersects $\varphi_{\mathbf{i}}(I)$.

The point ii) of Definition \ref{def:conda} can then be written as

\[
\forall_{x\in I_2^\alpha} \ \sum_{\mathbf{i} \in D_{r_\alpha}(x, I^\alpha _1, \alpha)} p_{\mathbf{i}} \geq 2.
\]
\begin{figure}[ht!]
  % Requires \usepackage{graphicx}
  \includegraphics[width=12cm]{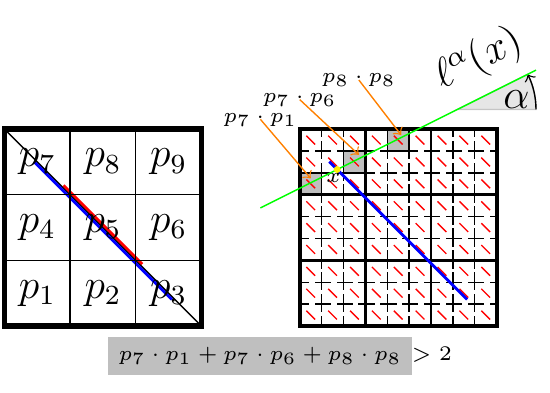}\\

 \caption{Condition $A(\alpha )$:  $I_{1}^{\alpha }$ is the small red,  $I_{2}^{\alpha }$ is the big blue interval on the left. $r_\alpha =2$ and
  the small red intervals are the scaled copies of $I_{1}^{\alpha }$.}\label{217}
\end{figure}

In other words, Condition $A(\alpha)$ is satisfied if for given $\alpha$ one can define 'small central' and 'large central' part of each square in such a way that for some $r\in \mathbb{N}$ if a line in direction $\alpha$ intersects the 'large central' part of some $n$-th level square then the expected number of 'small central' parts of its $n+r$-th level subsquares it intersects is uniformly greater than 1.

\subsection{Consequences of Condition $A(\alpha)$}

It is clear that if $A(\alpha)$ holds then one can apply the large deviation argument for projection in direction $\alpha$ - modulo a minor technical problem that the random variables in the large deviations theorem are not identically distributed.

A bit more complicated is the proof that almost surely all the projections contain intervals. It is based on the following robustness properties:

\begin{proposition} \label{prop:arob}
If condition A($\alpha$) holds for some $\alpha\in\mathcal{D}$ for some $I_1^\alpha, I_2^\alpha$ and $r_\alpha$ then it will also hold in some neighbourhood $J\ni\alpha$. Moreover, for all $\theta\in J$ we can choose $I_1^\theta=I_1', I_2^\theta=I_2, r_\theta=r_\alpha$ not depending on $\theta$.
\end{proposition}

A natural corollary is that the whole range $\mathcal{D}$ can be presented as a countable union of closed intervals $J_i=[\alpha_i^-, \alpha_i^+]$ such that for each $i$ Condition A($\alpha$) holds for all $\alpha\in J_i$ with the same $I_1^i, I_2^i, r_i$.

\begin{proposition}
Let $I\subset B(I,\ell) \subset J \subset \Delta_\alpha$. If $\mathbf{i} \in D_n(x, I, \alpha)$ then $\mathbf{i} \in D_n(x, J, \beta)$ for all $\beta \in (\alpha - \ell M^{-n}, \alpha + \ell M^{-n})$.
\end{proposition}

Hence, inside each $J_i$ one does not need to repeat the large deviation argument separately for each $\alpha$. At level $n$ it is enough to check it for approximately $M^n$ directions. As the number of directions one needs to check grows only exponentially fast with $n$, the proof goes through.

\subsection{Checking Condition $A(\alpha)$}

One last thing needed is an efficient way to check whether $A(\alpha)$ holds.

\begin{definition}
We say that the fractal percolation model satisfies {\bf Condition B($\alpha$)} if there exists a nonnegative continuous function $f:\Delta_\alpha\to \R$ such that $f$ is strictly positive except at the endpoints of $\Delta_\alpha$ and that
\begin{equation} \label{25}
F_\alpha f \geq (1+\varepsilon)f
\end{equation}
for some $\varepsilon>0$.
\end{definition}

\begin{proposition} \label{prop:ba}
B($\alpha$) implies A($\alpha$).
\end{proposition}

In particular, for homogeneous case $p_i=p > M^{-1}$ for any $\alpha$ one can choose $f(x)$ as the length of the intersection of $Q$ with the line in direction $\alpha$ passing through $x$. It is easy to check that this function satisfies \eqref{25} for $\varepsilon = pM-1$.

\subsection{Application: visibility}

For a given set $E$, we define the visible subset (from direction $\alpha$) as the set of points $x\in E$ such that the half-line starting at $x$ and going in direction $\alpha$ does not meet any other point $y\in E$. Similarly, given $z\in \mathbb{R}^2$, the visible subset (from $z$) is the set of points $x\in E$ such that the interval $\overline{xz}$ does not meet any other point $y\in E$.

Let $E$ be a homogeneous fractal percolation with $p>M^{-1}$. By Theorem \ref{thm:ortho}, $E$ is quite opaque: the orthogonal projection in any direction almost surely contain intervals. In particular, with large probability it contains large intervals. By stochastical self-similarity of $E$, the same is true for each $E\cap Q_{\mathbf{i}}$. Hence, not many points can be visible:

\begin{theorem}[ \cite{Arhosalo2012}]
If $E$ is nonempty, almost surely the visible set from direction $\alpha$ has finite one-dimensional Hausdorff measure for each $\alpha$ and the visible set from point $z$ has Hausdorff dimension 1 for each $z\in \mathbb{R}^2$.
\end{theorem}

\section{General projections: the transparent case}

%The results reviewed in the previous section imply that whenever the homogeneous fractal percolation set $E$ satisfies $\dim_{\rm H}E>1$ (that is, $p>M^{-1}$), then  for every $\alpha $, $\mathrm{proj}_\alpha E$ is as big as it can be, that is it contain some intervals.
%In this section we assert that similarly, in the case when $\dim_{\rm H}E\leq 1$ (that is $M^{-2}<p\leq M^{-1}$) the projection of the homogeneous fractal percolation set is as big as it can be. More precisely,

In this section we present results analogous to the second part of the Marstrand theorem. For homogeneous fractal percolation with Hausdorff dimension smaller than 1 almost surely $\dim_{\rm H}(\mathrm{proj}_\alpha(E)) =\dim_{\rm H} E$ for all $\alpha$. Together with the results of the previous section, it implies

\begin{theorem}[\cite{RamsII}]\label{n5}
In the homogeneous case, that is $E=E^h(2,M,p)$ for almost all realizations of $E$
\begin{equation}\label{u2}
\forall \alpha ,\  \dim_{\rm H}(\mathrm{proj}_\alpha E)=\min\left\{1,\dim_{\rm H}(E)\right\}.
\end{equation}
\end{theorem}
\textbf{Principal Assumption for this Section:} In this section we always work in the homogeneous case:
$$
E=E^h(2,M,p),
$$
where
\begin{equation}\label{u3}
 M^{-2}<p\leq M^{-1}.
\end{equation}
That is $p$ is chosen to ensure that $E\ne\emptyset $ with positive probability and
$\dim_{\rm H}(E)\leq 1$ almost surely conditioned on non-extinction.
To prove Theorem \ref{n5} one needs to analyze the structure of the slices of $E_n$:

\medskip

\textbf{Informal description of the structure of slices of} $\mathbf{E_n}$ (which was defined as the $n$-th approximation of $E$):
Namely,
    for almost all realizations of $E $  and for \textbf{all} straight lines $\ell$:
   the number of level $n$ squares having nonempty intersection with $E$ is at most $\mathrm{const}\cdot n$.
On the other hand, almost surely for $n$ big enough, we can find some line of $45^{\circ}$ angle which intersects $\mathrm{const}\cdot n$ level $n$ squares.

\medskip

Let $\mathcal{L}^\varepsilon $ be the set of lines on the plane whose angle is separated both from $0^{\circ}$ and $90^{\circ}$ at least by $\varepsilon $. Further for a line $\ell$ let $\mathcal{E}_n(\ell)$ be the set of retained level $n$ squares that intersect $\ell$. That is,
$$
\mathcal{E}_n(\ell):=\left\{\mathbf{i}\in \mathcal{E}_n:Q_{\mathbf{i}}\cap \ell\ne\emptyset \right\}.
$$

\begin{theorem}[\cite{RamsII}]\label{u66}
For almost all realizations of $E$  we have
\begin{equation}\label{v6}
\forall \varepsilon  \in \left(0,\frac{\pi }{2}\right),\
\exists N,\ \forall n\geq N,\ \forall \ell\in \mathfrak{L}^\varepsilon ; \quad
\#\mathcal{E}_n(\ell)\leq \mathrm{const}\cdot n.
\end{equation}

\end{theorem}

For simplicity, the proof in horizontal/vertical direction only (for general directions one needs to apply techniques presented in previous subsection). The proof is once again based on the large deviation argument, but working in the opposite direction. This time the expected number of squares in a subcolumn is smaller (by a constant bounded away from 1) than the number of squares in the column (and not greater, like in the opaque case). Hence, we can guarantee that if the column has sufficiently many squares for the large deviation theory to work, the number of squares in all subcolumns will shrink. This leads to an estimation on the possible rate of growth.

This estimation is sharp:

\begin{proposition}[\cite{RamsII}]\label{u76}
There exists a constant $0<\lambda <1$ such that
for almost all realizations, conditioned on $E\ne\emptyset $, there exists an $N$ such that for all
$n>N$ there exists a line $\ell$ with
\begin{equation}\label{u77}
  \#\mathcal{E}_n\left(\ell\right)>\lambda n.
\end{equation}
\end{proposition}

Theorem \ref{n5} is an immediate consequence of Theorem \ref{u66}.

\section{The arithmetic sum of at least three fractal percolations}

To study arithmetic sums of more than two fractal percolations we need to combine results of the previous three sections. Like in section 4, we look at the projection $(x_1,\ldots,x_d) \to \sum x_i$ from the cartesian product of fractal percolations to the real line. The proof is based on the large deviation argument presented in section 5. However, the main technical difficulty is the presence of dependencies. We will use the results from section 6 to bound their impact.

%In this section we explain how to use a result similar to Theorem \ref{u66} to prove the existence of some intervals in the arithmetic sum of three homogeneous fractal percolations.
Let
$$
E^{i}:=E^h\left(1,M,p_i\right),\ i=1,2,3,\quad
p:=p_1\cdot p_2\cdot p_3\mbox{ and }
E:=E^h\left(3,M,p\right).
$$
Then
$$
\dim_{\rm H}\left(E^1\times E^2\times E^3\right)=\dim_{\rm H}(E)=\frac{\log M^3\cdot p}{\log M}.
$$
  Moreover, the probability that a level $n$ cube $C$ is contained in any of the two random Cantor sets above is equal to $p^n$.

Let $S_a$ be the plane $\{\sum x_i =a\}$. We can write
 $$
 E^{\mathrm{sum}}:=
E^1+ E^2+E^3=
\left\{a:S_a\cap \left(E^1\times E^2\times E^3\right)\ne\emptyset \right\}.
$$

 \begin{wrapfigure}{r}{0.4\textwidth}
  \begin{center}
    \includegraphics[width=5cm]{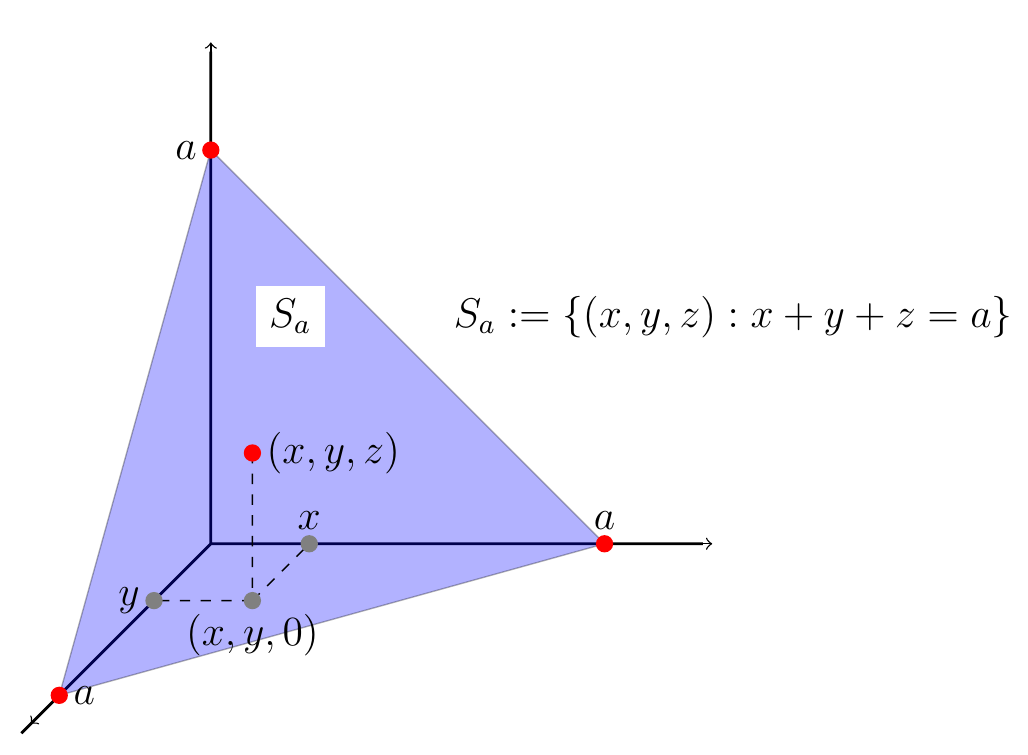}
  \end{center}
 % \caption{fan}
\end{wrapfigure}
 That is we can consider  $ E^{\mathrm{sum}} $ as the projection of $E_1\times E _2\times E _3 $ to the $x$-axis with planes orthogonal the vector $(1,1,1)$. So,
 $ E^{\mathrm{sum}} $ can contain an interval only if its dimension is greater than one, that is $p>M^{-2}$.
 %In this case however it does contain an interval. We proved actually much more:
It is a sufficient condition as well:
\begin{theorem}[\cite{RamsII}] \label{thm:prod} Let $d\geq 2$ and for $i=1,\dots ,d$ let $E^{i}:=E^h(1,M,p_i)$ satisfying
\begin{equation} \label{cond1}
p:=\prod_{i=1}^{d} p_i > M^{-d+1}.
\end{equation}
Then for every $\mathbf{b}=(b_1,\dots ,b_d)\in \mathbb{R}^d$, $b_i\ne 0$ for all $i=1,\dots ,d$
 the sum $ E^{\mathrm{sum}} _\mathbf{b}=\sum\limits_{i=1}^{d} b_i {E}^{i}$ contains an interval almost surely, conditioned on all ${E}^{i}$ being nonempty.
\end{theorem}
We explain the proof of this theorem in the special case when $d=3$ and $\mathbf{b}=(1,1,1)$. To verify that a certain $a\in  E^{\mathrm{sum}} $ we need to prove that the $n$ approximation of the product intersects $S_a$, that is $(E^1\times E^2\times E^3)_{n}\cap S_a\ne\emptyset $ for every $n$.
 It follows from the dimension formula and \eqref{cond1}
that we have  $M^{n(1+\tau )}$ retained level $n$ cubes for some $\tau >0$. By the pigeon hole principle for at least one $k=0,\dots ,3M^n$ the plane $S_{kM^{-n}}$  intersects at least $M^{n\tau }$ retained level $n$ cubes. For such a $k$ we write $a=kM^{-n}$. So, $\#\left\{\mathcal{E}_n\cap S_a\right\}\geq M^{n\tau }$.

Fix an $0\leq m\leq M $. How many level $n+1$ retained cubes  intersect
 $S_{a+mM^{-(n+1)}}$? If the way $E^1\times E^2\times E^3$ develops in every level $n$ cube was independent then we could get that the answer by the large deviation argument: exponentially many except for an event with a super exponentially small probability.

\begin{figure} \label{f7}
  % Requires \usepackage{graphicx}
  \includegraphics[height=5cm]{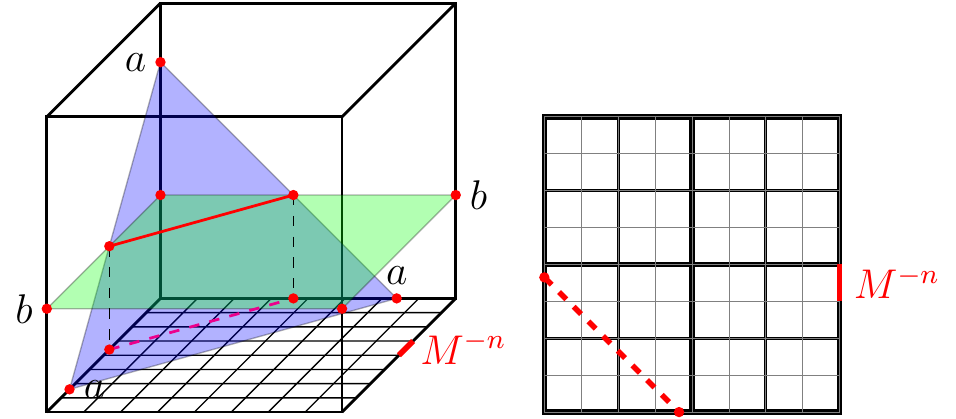}\\
  \caption{The cubes intersecting the red line are not independent}\label{n8}
\end{figure}

We remind that the cubes are dependent if they have the same $x_1, x_2$ or $x_3$ coordinate. Figure \ref{f7} shows the geometric position of (some of: we consider only the cubes with the same $x_3$ coordinate) cubes dependent on one chosen cube: $x_1+x_2+x_3={\rm const}$ and $x_3={\rm const}$ imply $x_1+x_2={\rm const}$. Potentially there could be exponentially many
such cubes. The key step of the proof is that
using a theorem analogous to Theorem \ref{u66} for $E^1\times E^2$ instead of
$E^h(2,M,p_1\cdot p_2)$ one can check that on the red dashed line on Figure
\ref{n8} there are only constant times $n$ retained squares, consequently
the $M^{n\tau }$ level $n$ cubes having non-empty intersection with $S_{a}$ (the blue
plane on Figure \ref{n8}) can be divided into $\mathrm{const}\cdot n$ classes such that the coordinate axes projection of any two cubes in a class are different. The events inside each class are independent, hence we can use the large deviation theory separately for each class.
A technical comment: in order to be able to go with this procedure we may have to decrease $p_1,p_2,p_3$ in such a way that for the modified values we have
$$
p_1\cdot p_2\cdot p_3>M^{-2}\mbox{ but } p_i\cdot p_j<M^{-1} \mbox{ for distinct } i,j\in \left\{1,2,3\right\}.
$$
That is, $E^1\times E^2\times E^3$ is a big set in the sense that it has dimension greater than one but its all coordinate plane projections should be small sets having dimension smaller than one -- only then the $n$-th approximates of the coordinate plane projections intersect every line in at most $\mathrm{const}\cdot n$ retained squares. However, the property of almost surely having intervals in the algebraic sum is monotonous with respect to $\{p_i\}$.

Hence among those level $n$ retained cubes that intersect the blue plane $S_a$ there cannot be more than $\mathrm{const}\cdot n$ on the red line (any coordinate plane parallel line) which imply that the number of cubes dependent on any one cube is polynomial ($\mathrm{const}\cdot n$). This bound on the dependency matrix lets us control the dependencies.

%and this gives a way to fix the problem caused by the dependence between the  portions of $E^1\times E^2\times E^3$ in retained level $n$ cubes having common coordinate axes projections.

\bibliographystyle{plain}

\bibliography{biblo_5}

\begin{thebibliography}{10}

\bibitem{Arhosalo2012}
I.~Arhosalo, E.~J{\"a}rvenp{\"a}{\"a}, M.~J{\"a}rvenp{\"a}{\"a}, M.~Rams, and
  P.~Shmerkin.
\newblock Visible parts of fractal percolation.
\newblock {\em Proceedings of the Edinburgh Mathematical Society (Series 2)},
  55(02):311--331, 2012.

\bibitem{Chayes1988}
J.T. Chayes, L.~Chayes, and R.~Durrett.
\newblock Connectivity properties of mandelbrot's percolation process.
\newblock {\em Probability theory and related fields}, 77(3):307--324, 1988.

\bibitem{Chayes1996}
L.~Chayes.
\newblock On the length of the shortest crossing in the super-critical phase of
  mandelbrot's percolation process.
\newblock {\em Stochastic processes and their applications}, 61(1):25--43,
  1996.

\bibitem{Dekking1988}
F.M. Dekking and G.R. Grimmett.
\newblock Superbranching processes and projections of random cantor sets.
\newblock {\em Probability theory and related fields}, 78(3):335--355, 1988.

\bibitem{Dekking2009}
M.~Dekking.
\newblock Random cantor sets and their projections.
\newblock {\em Fractal Geometry and Stochastics IV}, pages 269--284, 2009.

\bibitem{Dekking1990}
M.~Dekking and R.W.J. Meester.
\newblock On the structure of mandelbrot's percolation process and other random
  cantor sets.
\newblock {\em Journal of Statistical Physics}, 58(5):1109--1126, 1990.

\bibitem{Dekking2008}
M.~Dekking and K.~Simon.
\newblock On the size of the algebraic difference of two random cantor sets.
\newblock {\em Random Structures \& Algorithms}, 32(2):205--222, 2008.

\bibitem{Falconer1986}
K.J. Falconer.
\newblock Random fractals.
\newblock {\em Math. Proc. Cambridge Philos. Soc}, 100(3):559--582, 1986.

\bibitem{Falconer1992}
K.J. Falconer and G.R. Grimmett.
\newblock On the geometry of random cantor sets and fractal percolation.
\newblock {\em Journal of Theoretical Probability}, 5(3):465--485, 1992.

\bibitem{Falconer1994}
K.J. Falconer and G.R. Grimmett.
\newblock Correction: On the geometry of random cantor sets and fractal
  percolation.
\newblock {\em Journal of Theoretical Probability}, 7(1):209--210, 1994.

\bibitem{Hawkes1981}
J.~Hawkes.
\newblock Trees generated by a simple branching process.
\newblock {\em Journal of the London Mathematical Society}, 2(2):373--384,
  1981.

\bibitem{Kahane1976}
J.-P. Kahane and J.~Peyriere.
\newblock Sur certaines martingales de benoit mandelbrot.
\newblock {\em Advances in mathematics}, 22(2):131--145, 1976.

\bibitem{Mandelbrot1974}
B.B. Mandelbrot.
\newblock Intermittent turbulence in self-similar cascades- divergence of high
  moments and dimension of the carrier.
\newblock {\em Journal of Fluid Mechanics}, 62(2):331--358, 1974.

\bibitem{Mandelbrot1983}
B.B. Mandelbrot.
\newblock The fractal geometry of nature/revised and enlarged edition.
\newblock {\em New York, WH Freeman and Co., 1983, 495 p.}, 1, 1983.

\bibitem{Marstrand1954}
J.~M. Marstrand.
\newblock Some fundamental geometrical properties of plane sets of fractional
  dimensions.
\newblock {\em Proceedings of the London Mathematical Society}, 3(1):257--302,
  1954.

\bibitem{Mauldin1986}
R.D. Mauldin and S.C. Williams.
\newblock Random recursive constructions: asymptotic geometric and topological
  properties.
\newblock {\em Trans. Amer. Math. Soc}, 295(1):325--346, 1986.

\bibitem{Mora2009}
P.~Mora, K.~Simon, and B.~Solomyak.
\newblock The lebesgue measure of the algebraic difference of two random cantor
  sets.
\newblock {\em Indagationes Mathematicae}, 20(1):131--149, 2009.

\bibitem{RamsII}
M.~Rams and K.~Simon.
\newblock The dimension of projections of fractal percolations.
\newblock {\em preprint}.

\bibitem{Rams}
M.~Rams and K.~Simon.
\newblock Projections of fractal percolations.
\newblock {\em To appear in Ergodic Theory and Dynamical Systems}.

\end{thebibliography}

\end{document}